\documentclass[12pt,longbibliography]{article}

\usepackage[utf8]{inputenc}

\usepackage[T1]{fontenc}

\usepackage[a4paper, margin=2.7cm]{geometry}

\usepackage{amsmath, amsthm, amsfonts, amssymb}
\usepackage{mathrsfs}   
\usepackage{textcomp}        

\usepackage[hyphens]{url} 
\usepackage[all]{xy}

\usepackage{amscd}
\usepackage[colorlinks=true,linkcolor=blue,citecolor=blue,urlcolor=blue,breaklinks]{hyperref}

\newtheorem{innercustomgeneric}{\customgenericname}
\providecommand{\customgenericname}{}
\newcommand{\newcustomtheorem}[2]{%
	\newenvironment{#1}[1]
	{%
		\renewcommand\customgenericname{#2}%
		\renewcommand\theinnercustomgeneric{##1}%
		\innercustomgeneric
	}
	{\endinnercustomgeneric}
}

\newcustomtheorem{customthm}{Theorem}

\theoremstyle{plain}
\newtheorem{theorem}{Theorem}[section]
\newtheorem{lem}[theorem]{Lemma}
\newtheorem{cor}[theorem]{Corollary}
\newtheorem{prop}[theorem]{Proposition}

\theoremstyle{definition}
\newtheorem{defn}[theorem]{Definition}
\newtheorem{example}[theorem]{Example}

\theoremstyle{remark}
\newtheorem{remark}[theorem]{Remark}

\newcommand{\R}{\mathbb{R}}
\newcommand{\N}{\mathbb{N}}

\newcommand{\Q}{\mathbb{Q}}
\newcommand{\Z}{\mathbb{Z}}

\hypersetup{pdftitle={Title of the PDF}}
\hypersetup{pdfauthor={Author of the PDF}}

\title{Profinite genus of fundamental groups
	of torus bundles}
\author{Genildo de Jesus Nery\thanks{The author held CAPES scholarship during the preparation of this article.}}

\AtEndDocument{\bigskip{\footnotesize
		\noindent
		\textsc{Genildo de Jesus Nery. University of Brasilia, Department of Mathematics, Brasilia DF, Brazil.}
		\textit{E-mail address}: \texttt{\href{mailto: genildo.j.nery@gmail.com}{genildo.j.nery@gmail.com}}
}}

\begin{document}
	
	\maketitle
	
	\begin{abstract}
		In this paper we establish lower and upper bounds for the cardinality of the profinite genus of the fundamental group $\pi_{1}(M_A)\cong (\Z \times \Z)\rtimes_{A}\Z$ of a torus bundle $M_{A}$ in terms of the number of ideal classes of the order $\Z[\lambda]$, where $\lambda$ is an eigenvalue of the matrix $A$ in $\mathrm{GL}_{2}(\Z)$. \\
		\
		\
	\end{abstract}

	
	\section{Introduction}
	\label{sec:intro}
	
	Recently, the following question has attracted much
	attention in geometric group theory.
	
	\begin{flushleft} 
		\textbf{Question 1:} Let $M$ be a compact, orientable 3-manifold. To what extent is $\pi_{1}(M)$ determined by its profinite completion?
	\end{flushleft}
	In particular, \cite[p. 138]{AFW15} highlights the following conjecture.
	
	\begin{flushleft}
		\textbf{Conjecture:} The fundamental group of a compact, orientable 3-manifold is determined by its profinite completion up to finitely many isomorphism classes.
	\end{flushleft}

	W. P. Thurston, in 1982, has conjectured that there are only eight 3-dimensional model geometries, which are: $\mathbf{E}^{3}$,  $\mathbf{H}^{3}$, $\mathbf{S}^{3}$, $\mathbf{S}^{2}\times \R$, $\mathbf{H}^{2}\times\R$, $\widetilde{\mathbf{SL}_{2}}(\R)$, $\mathbf{Nil}$ and $\mathbf{Sol}$; this was proved by  G. Y. Perelman in 2003. H. Wilton and P. Zalesskii \cite{WZ17} showed that the geometry of a  $3$-manifold is determined by the profinite completion of its fundamental group. In \cite{WZ18}, they also proved that the profinite completion of the fundamental group of a 3-manifold $M$ determines the Jaco--Shalen--Johannson decomposition of $M$.  H. Wilton and P. Zalesskii highlight that considering \cite{WZ17} and \cite{WZ18}, the next step in addressing Question 1 is to consider the parts of the Jaco--Shalen--Johannson decomposition (i.e., the manifold modeled in one of the eight geometries) and point out that a definitive treatment of the case of 3-dimensional solvmanifolds would be a valuable addition to the literature. 
	
	G. Wilkes \cite{Wil17} proved that all Seifert fiber spaces are distinguished from each other by their profinite completions apart from some previously-known examples due to Hempel of 3-manifolds with the geometry $\mathbf{H}^{2}\times\R$. Furthermore, in \cite{Wil18}, he presented a complete answer to the Question 1 in relation to graph manifolds, not including the 3-dimensional solvmanifolds. In other words, G. Wilkes proved that a closed orientable 3-manifold with geometric structure modeled in one of the following geometries: $\mathbf{S}^{3}$, $\mathbf{E}^{3}$, $\mathbf{S}^{2}\times \R$, $\widetilde{\mathbf{SL}_{2}}(\R)$ or $\mathbf{Nil}$, is determined by the profinite completion of its fundamental group. By contrast, L. Funar \cite{Fun12} constructed a family of torus bundles, with the geometry $\mathbf{Sol}$, which are not distinguished by the profinite completions of their fundamental groups. 
	
	Here we give an  answer to Question 1 for $3$-dimensional solvmanifolds addressing the remark of H. Wilton and P. Zalesskii \cite{WZ18}. We prove that there is a family of torus bundles, with the geometry $\mathbf{Sol}$, that  are determined by the profinite completions of their fundamental groups.

	A solvmanifold can be represented by a Mostow fibration (see for instance \cite[Thm. 2.4]{JKM08}). For a 3-dimensional solvmanifold $M$ that is not a nilmanifold, the corresponding
	Mostow fibration must have one of the following two forms,
	\begin{eqnarray*}
		\mathbb{T}^{2} \hookrightarrow M\rightarrow S^{1} \ \ \text{or} \ \ S^{1}\hookrightarrow M \rightarrow \mathbb{T}^{2}, 
	\end{eqnarray*}
	where $S^{1}$ is the unit circle and $\mathbb{T}^{2}$ is a 2-dimensional torus. The $3$-dimensional solvmanifolds $M$ of type $S^{1}\hookrightarrow M \rightarrow \mathbb{T}^{2}$ are determined by the profinite completions of their fundamental groups (see  \cite[Thm. 1.2]{Wil17}). Thus we concentrate on the case $\mathbb{T}^{2} \hookrightarrow M\rightarrow S^{1}$ that leads to corresponding exact sequence of the fundamental groups, 
	\begin{eqnarray*}
		0\rightarrow \Z^{2} \rightarrow\pi_{1}(M)\rightarrow \Z \rightarrow 0. 
	\end{eqnarray*}
	
	A $3$-dimensional solvmanifold $M$  of type $\mathbb{T}^{2} \hookrightarrow M\rightarrow S^{1}$ is called a \textit{torus bundle}. Its fundamental group is the semi-direct product $\Z^{2}\rtimes_{A}\Z$, where $A$ is a matrix in $\mathrm{GL}_{2}(\Z)$. It is well known that a compact solvmanifold is uniquely determined by its fundamental group, up to diffeomorphism (see \cite[Thm. A]{Mos54}). Thus, it makes sense to write $M_{A}$ to denote the torus bundle $M$ of fundamental group $\Z^{2}\rtimes_{A}\Z$, with $A$ in $\mathrm{GL}_{2}(\Z)$.
	
	We will denote by $\mathrm{tr}(A)$ and $\mathrm{det}(A)$ the trace and determinant of a matrix $A$. If $|\mathrm{tr}(A)|> 2$, then $A$ is hyperbolic (i.e., neither of its eigenvalues has absolute value $1$) and the torus bundle $M_{A}$ admits the $\mathbf{Sol}$ geometry; if $|\mathrm{tr}(A)|\leq 2$,  the torus bundle $M_A$ admits the $\mathbf{Nil}$ geometry (respectively the $\mathbf{Sol}$ geometry) if $A$ has infinite order  and is not hyperbolic (respectively $A$ is a hyperbolic matrix); and if $A$ has finite order, the torus bundle $M_A$ admits the $\mathbf{E}^{3}$ geometry (see \cite[Thm. 5.5]{Sco83}).

	We denote by $\widehat{G}$ the profinite completion of a group $G$. Let $\mathfrak{C}$ be a class of groups and $G\in \mathfrak{C}$. Following \cite{GZ11}, we define the $\mathfrak{C}$\textit{-genus} of $G$, denoted by $\mathfrak{g}(\mathfrak{C},G)$, as the set of isomorphism classes of groups belonging to $\mathfrak{C}$ that has the same profinite completion as the fixed group $G$. Since the fundamental groups of the torus bundles are polycyclic, in this paper we restrict attention to the class $\mathcal{PF}$ of  polycyclic-by-finite groups. For this reason, we calculate the profinite genus $\mathfrak{g}(\mathcal{PF},G_{A})$  for any fixed group $G_{A}=\Z^{2}\rtimes_{A} \Z$, where $A\in \mathrm{GL}_{2}(\Z)$. In what follows, $\#\mathfrak{g}(\mathfrak{C},G)$ denotes the cardinality of the set $\mathfrak{g}(\mathfrak{C},G)$.
	
	In our study of the Question 1, we start with the case where $\Z^{2}\rtimes_{A} \Z$ is nilpotent.
	
	\begin{theorem}\label{Caso nilpotente Th}
		Let $A$ be a matrix in $\mathrm{GL}_{2}(\Z)$ and consider the semi-direct product $G_A=\Z^{2}\rtimes_{A} \Z$. If all eigenvalues of $A$ are equal to $1$, then $\#\mathfrak{g}(\mathcal{PF},G_A)=1$.
	\end{theorem}
	
	In the case that $\Z^{2}\rtimes_{A} \Z$ is not nilpotent, we establish lower and upper bounds for the cardinality of $\mathfrak{g}(\mathcal{PF},\pi_1(M_A))$  in terms of the number of ideal classes of the order $\Z[\lambda]$, where $\lambda$ is an eigenvalue of the matrix $A$ in $\mathrm{ GL}_{2}(\Z)$. To this end, we use a famous theorem of Latimer and MacDuffee \cite[Thm. III. 13]{New72}, which says that there is a one-to-one correspondence between the conjugacy classes of $n\times n$ matrices over $\Z$ with same characteristic polynomial $f$ and the ideal classes of the order $\Z[\lambda]$, where $\lambda$ is a root of $f$. More precisely, we get the following result.
	
	\begin{theorem}\label{Caso nao nilpotente Th}
		Let $A$ be a matrix in $\mathrm{GL}_{2}(\Z)$ and consider the group $G_{A}=\Z^{2}\rtimes_{A}\Z$.  Let $\lambda$ be an eigenvalue of $A$. Then,
		
		\begin{enumerate}
			\item[(i)] if $A$ has distinct eigenvalues, $\mathrm{tr}(A)\neq 0$ and the conjugacy class of $A$ corresponds to a class of an invertible ideal of the order $\Z[\lambda]$, then
			
			\begin{equation*}
				\left\{
				\begin{array}{c c}
					h(\Z[\lambda])\leq \#\mathfrak{g}(\mathcal{PF}, G_A)\leq \tilde{h}(\lambda), & \text{if} \ \ \mathrm{det}(A)=-1\\
					h(\Z[\lambda])/2 \leq \#\mathfrak{g}(\mathcal{PF},G_A)\leq \tilde{h}(\lambda), & \text{if} \  \ \mathrm{det}(A)=1	
				\end{array}
				\right.
			\end{equation*}
			where $h(\Z[\lambda])$ is the order of the ideal class group and $\tilde{h}(\lambda)$ is the number of ideal classes of $\Z[\lambda]$.
			
			\item[(ii)] if $\mathrm{tr}(A)=0$ or $A$ has all its eigenvalues equal to $-1$, then $$\#\mathfrak{g}(\mathcal{PF}, G_{A})= h(\Z[\lambda])=1.$$
		\end{enumerate}
	\end{theorem}
	
	Let us mention two important consequences of the Theorems \ref{Caso nilpotente Th} and \ref{Caso nao nilpotente Th}. 
	
	\begin{cor}\label{cor1}
		Let $A\in \mathrm{GL}_{2}(\Z)$ with $\mathrm{det}(A)=-1$ and consider the semi-direct product $G_{A}=\Z^{2}\rtimes_ {A} \Z$. If $\mathrm{tr}(A)^{2}-4\mathrm{det}(A)$ is square-free, then $\#\mathfrak{g}(\mathcal{PF}, G_{A})= h(\Z[\lambda])$ where $\lambda$ is an eigenvalue of $A$.
	\end{cor}
	
	\begin{cor}\label{cor2}
		Let $A$ be a matrix in $\mathrm{GL}_{2}(\Z)$, as in Theorem \ref{Caso nao nilpotente Th} or Theorem \ref{Caso nilpotente Th}, with eigenvalue $\lambda$. If $\tilde{h}(\lambda)=1$, then the torus bundle $M_A$ is determined among  $3$-manifolds by the profinite completion of its fundamental group.	
	\end{cor}  
	
	It is worth pointing out that, in 1801, Gauss \cite{Gau66} conjectured --and this remains open to this day-- that there are infinitely many real quadratic fields with class number one. 
	
	The Theorems \ref{Caso nilpotente Th} and \ref{Caso nao nilpotente Th} are proved in Section \ref{main results}. Sections \ref{Preliminary Results}, \ref{The Ideal Class Group} and \ref{Conjugation in GL} contain preliminary results and definitions necessary for the proof of the main theorems.  We conclude the paper in Section \ref{examples} with examples of applications of the Theorems  \ref{Caso nilpotente Th} and \ref{Caso nao nilpotente Th}.
	
	\section{Preliminary Results}\label{Preliminary Results}
	
	Let $A$ be  a matrix in $\mathrm{GL}_2(\Z)$. Consider the semi-direct product $G_A=N\rtimes_{A} C$, where  $N\cong \Z\times \Z$ and $C \cong \Z$ so that the generator $1\in C$ acts on $N$ via $A$.
	
	\begin{lem}\label{(id-A)(N)=[G,G]}
		Let $A$ be a matrix in $\mathrm{GL}_{2}(\Z)$ such that none of its eigenvalues  is $1$. Let $G_A=N\rtimes_{A} C$ be the corresponding semi-direct product. Then,
		\begin{enumerate}
			\item[(i)] $(\mathrm{id}-A)(N)$ is a subgroup of finite index in $N$.
			\item[(ii)] The derived subgroup of $G_{A}$ is equal to $(\mathrm{id}-A)(N)$.
		\end{enumerate} 
	\end{lem}
	\begin{proof}
		We write the elements of $G_{A}$ as pairs $(v,t)\in N\times C$ with multiplication
		\begin{equation*}
			(v,t)(u,s)=(v+A^{t}u,t+s).
		\end{equation*}	
		\begin{enumerate}
			\item[(i)] Consider the mapping $f:N\rightarrow N$ given by $f(v)= v-Av$. It is easy to see that $f$ is a homomorphism. Therefore, $f(N)=(\mathrm{id}-A)(N)$ is a subgroup of $N$. Since none of the eigenvalues of the matrix $A$  is equal to $1$, the equation $Av=v$ has only the trivial solution $v=0$. This implies that $f$ is injective. Hence, $(\mathrm{id}-A)(N)\cong n\Z \times m\Z$ for some nonzero integer $n$ and $m$, and so $(\mathrm{id}-A)(N)$ has finite index in $N$.
			
			\item[(ii)] Note that the quotient
			\begin{equation*}
				\frac{G_A}{(\mathrm{id}-A)(N)}=\frac{N}{(\mathrm{id}-A)(N)}\rtimes_{A} C
			\end{equation*}
			is an abelian group, since $A$ acts trivially on $N/(\mathrm{id}-A)(N)$. Thus, $$[G_{A},G_{A}]\leq (\mathrm{id}-A)(N).$$ On the other hand, we have 
			\begin{eqnarray*}
				(v-Av,0)&=& (v,t)(Av,t)^{-1} \\
				&=& (v,t)\left((0,1)(v,t)(0,1)^{-1} \right)^{-1} \\
				&=& (v,t)(0,1)(v,t)^{-1}(0,1)^{-1},
			\end{eqnarray*}
			for each $v\in N$ and $t\in C$. Then, $(\mathrm{id}-A)(N)\subseteq [G_A,G_A]$ and therefore $[G_A,G_A]=(\mathrm{id}-A)(N)$.
		\end{enumerate}
	\end{proof}
	
	\begin{lem}\label{G_A is nipotent iff A has all the eigenvalues equal to 1}
		Let $A\neq I$ be a matrix in $\mathrm{GL}_{2}(\Z)$. Then, $G_A$ is a nilpotent group if and only if $A$ has all the eigenvalues equal to $1$. Furthermore, $G_A$ has  nilpotency class $2$.
	\end{lem}
	\begin{proof}
		Note that, if the matrix $A$ has at least one eigenvalue different from $1$, then the center of the group $G_A$ is trivial and, in consequence, $G_A$ is not nilpotent.
		
		Conversely, let $A\in \mathrm{GL}_{2}(\Z)$ with all its eigenvalues equal to $1$. By \cite[Lem. 13.27]{DK18},
		\begin{equation*}
			\{1\}\leq \Z \leq \Z^{2}
		\end{equation*}
		is a series of $\Z^{2}$ such that $A$ acts on $\Z^{2}/\Z$ and $\Z$ as the identity. Thus, $\Z$  is a central subgroup of $G_ {A}$, i.e., $\Z$ is contained in the center $Z_1(G_A)$ of $G_A$. Note that, $G_A$ is nilpotent of class $2$ if and only if we have the following upper central series
		\begin{equation*}
			\{1\}=Z_{0}(G_{A})\leq Z_1(G_A)\leq Z_2(G_A)=G_A,
		\end{equation*}
		if and only if the group $G_A/Z_1(G_A)$  is abelian. Which is the case, since the group
		\begin{equation*}
			\frac{\Z^{2}}{\Z}\rtimes_{A} \Z\cong \frac{G_A}{\Z}
		\end{equation*}
		is abelian (because $A$ acts trivially on $\Z^2/\Z$) and 
		\begin{equation*}
			\frac{G_A}{Z_1(G_A)}\cong \frac{G_A/\Z}{Z_1(G_A)/\Z}.
		\end{equation*}
	\end{proof}
	
	\begin{lem}\label{G_{A}=G_A^{-1}}
		Let $A$ be a matrix in $\mathrm{GL}_{2}(\Z)$. Then $G_{A}\cong G_{A^{-1}}$.
	\end{lem}
	\begin{proof}
		Consider the mapping $f:G_{A}\rightarrow G_{A^{-1}}$ defined by $f(v,t)=(v,-t)$. It is easy to see that $f$ is an isomorphism from $G_{A}$ to $G_{A^{-1}}$.
	\end{proof}
	
	Let $A$ and $B$ be matrices in $\mathrm{GL}_{2}(\Z)$. Consider the semi-direct products $G_{A}=N_{1}\rtimes_{A} C$ and $G_{B}=N_{2}\rtimes_{B} C$, where $N_{1}\cong \Z\times \Z \cong N_{2}$ and $C \cong \Z$.
	
	\begin{lem}\label{condicao2}
		\
		\begin{enumerate}
			\item[(i)] Let $A$ and $B$  be matrices in $\mathrm{GL}_{2}(\Z)$ such that none of its eigenvalues  is $1$. If $f:G_{A}\rightarrow G_{B}$ is an isomorphism then $f(N_{1})=N_{2}$.
			
			\item[(ii)] Let $A$ and $B$  be matrices in $\mathrm{GL}_{2}(\Z)$ such that none of its eigenvalues is $1$. If $f: \widehat{G}_{A}\rightarrow \widehat{G}_{B}$ is an isomorphism then $f(\widehat{N}_{1})= \widehat{N}_{2}$.
		\end{enumerate}
	\end{lem}
	\begin{proof} 
		\begin{enumerate}
			\item[(i)] We have by Lemma \ref{(id-A)(N)=[G,G]} that the derived subgroup $[G_B,G_B]$ has finite index in $N_2$. Consider
			\begin{equation*}
				L=\{\alpha\in G_B : \exists \ n\in \Z, n\neq 0, \ \text{such that} \ \alpha^{n}\in [G_B,G_B]\}.
			\end{equation*}
			We claim that $L=N_2$. Indeed, since the index of $[G_{B},G_{B}]$ in $N_{2}$ is finite we have  $N_{2}\subseteq L$. On the other hand, if $\alpha = (v, t) \in L \ (v\in N_{2}, t\in C)$, then there is $n\in \Z, \ n\neq 0$, such that $\alpha^{n}\in [G_{B},G_{B}]\leq N_{2}$. Note that, 
			\begin{equation*}
				\alpha^{n}=(v,t)^{n}=(v+A^{t}x+\cdots+A^{(n-1)t}v,nt).
			\end{equation*}
			Since $\alpha^{n}\in N_{2}$ and $n\neq 0$ we have $t=0$. Therefore $L \subseteq N_{2}$ and hence $L = N_{2}$. 
			
			Let $u\in N_{1}$, then there is $n\in \Z, \ n\neq 0$,  such that $u^{n}\in [G_{A},G_{A}]$. Hence, $f(u)^{n}\in [G_{B},G_{B}]$ and so $f(u)\in N_{2}$. Therefore, $f(N_{1})\subseteq N_{2}$. By a similar argument, $N_{2}\subseteq f(N_{1})$.
			
			\item[(ii)] This follows from (i).
		\end{enumerate}
	\end{proof}
	
	\begin{lem}[\cite{GZ11}, Prop 2.5]\label{prop 2.5}
		Let $H$, $N$ be groups and $\varphi_{1}, \varphi_{2}:H\rightarrow \mathrm{Aut}(N)$ be homomorphisms. Let $G_{1}=N\rtimes_{\varphi_{1}}H$ and $G_{2}=N\rtimes_{\varphi_{2}}H$ be the corresponding semi-direct products. The following statements hold.
		\begin{enumerate}
			\item Suppose that there exist $\Theta \in \mathrm{Aut}(H)$ and $\bar{f}\in\mathrm{Aut}(N)$ such that 
			\begin{eqnarray}\label{c1}
				\varphi_{1}(h)^{\bar{f}^{-1}}=\varphi_{2}(\Theta(h))
			\end{eqnarray}
			for all $h\in H$. Then the map $f: G_{1}\rightarrow G_{2}$ defined by $f(n,h)=(\bar{f}(n),\Theta(h))$ is an isomorphism which
			satisfies $f(N)=N$ and $f(H)=H$.
			\item Suppose now that $N$ is abelian and that $f: G_{1}\rightarrow G_{2}$ is an isomorphism such that $f(N)=N$. Define $\bar{f}=f_{|N}$ and $\Theta: H\rightarrow H$ by inducing $f$ to $G_{1}/N=H$. Then the pair $(\bar{f},\Theta)$ satisfies $\mathrm{(\ref{c1})}$ and hence also defines an isomorphism from $G_{1}$ to $G_{2}$.
		\end{enumerate}
	\end{lem}
	
	\begin{lem}[\cite{GZ11}, Cor 2.2]\label{cor 2.2}
		Let $H$, $N$ be groups and $\varphi_{1}, \varphi_{2}:H\rightarrow \mathrm{Aut}(N)$ be homomorphisms. Let $G_{1}=N\rtimes_{\varphi_{1}}H$ and $G_{2}=N\rtimes_{\varphi_{2}}H$ be the corresponding semi-direct products and let $f:G_{1}\rightarrow G_{2}$ be an isomorphism such that $f(N)= N$. Then $\widetilde{\varphi}_{1}(h)^{\widetilde{f}^{-1}}=\widetilde{\varphi}_{2}(f(h))$ for all $h\in H$, where $\widetilde{f}$ is the image of the restriction $f_{|N}$ in $\mathrm{Out}(N)$ and $\widetilde{\varphi}_{i}$ is the composition of $\varphi_{i}$ with the natural homomorphism $\mathrm{Aut}(N)\rightarrow \mathrm{Out}(N)$ for $i=1,2$. In particular, $\widetilde{\varphi}_{1}(H)$ and $\widetilde{\varphi_{2}}(H)$ are conjugate in $\mathrm{Out}(N)$.
	\end{lem}
	
	\begin{remark}
		Note that similar versions of Lemmas \ref{prop 2.5} and \ref{cor 2.2} hold for profinite groups and continuous homomorphisms. The same proofs presented in \cite{GZ11} applies to these reformulations.
	\end{remark}
	
	We will denote by $\mathfrak{F}(G)$ the set of finite quotients of a group $G$.
	
	\begin{prop}[\cite{FPDR82}, p. 227]
		Let $G$ and $H$ be finitely generated groups. Then $\widehat{G}$ and $\widehat{H}$ are isomorphic as topological groups if and only if $\mathfrak{F}(G)=\mathfrak{F}(H)$.
	\end{prop}
	
	This allows to state version  of Theorem of F. Grunewald and R. Scharlau.
	
	\begin{prop}[\cite{GS79}, p. 163]\label{Grunewald-Scharlau}
		Let $G_{1}$ and $G_{2}$ be finitely generated torsion-free nilpotent groups of nilpotency class 2, such that the Hirsch number \footnote{This is the number of infinite cyclic factors in a series with cyclic or finite factors.} of $G_{1}$ is smaller than or equal to 5 and  $\widehat{G}_{1}\cong\widehat{G}_{2}$. Then $G_{1}\cong G_{2}$.
	\end{prop}

	\section{The Ideal Class Group}\label{The Ideal Class Group}
	
	In this section we recall definitions and results of Number Theory that will be used in the paper.
	
	\begin{defn}
		Let $K$ be an algebraic number field of degree $n$. An \textit{order} of $K$ is a subring $\mathcal{O}$ of the ring of integers of $K$ which contains an integral basis of length $n$. The ring of integers of $K$, which we denote by $\mathcal{O}_{K}$, is called the maximal order of $K$.
	\end{defn}
	
	Set
	\begin{equation*}
		r:=\left\{
		\begin{array}{c c c}
			1, &\text{if} & m\not\equiv 1 \ (\mathrm{mod} \ 4)\\
			2, &\text{if} & m\equiv 1 \ (\mathrm{mod} \ 4)\\
		\end{array}
		\right.	
	\end{equation*}
	and define $\omega_0:=(r-1+\sqrt{m})/r$.
	
	\begin{prop}[\cite{JW09}, Thm. 4.17]\label{Thm 4.17 of JW09}
		If $\mathcal{O}$ is any order of a quadratic field $K$, then $\mathcal{O}$ is a $\Z$-module with basis $\{1,\omega\}$, where $\omega:=f\omega_0$ for some $f\in \Z$.
	\end{prop}
	
	Then $\mathcal{O}$ has finite index in $\mathcal{O}_K$. The index $f=|\mathcal{O}_K:\mathcal{O}|$ is called the \textit{conductor} of the order $\mathcal{O}$. 
	
	\newpage
	\begin{lem}\label{Existencia da unidade q nao e 1}
		Let $A$ be a matrix in $\mathrm{GL}_{2}(\Z)$. If $\mathrm{tr}(A)\neq 0$ and $A$ has distinct eigenvalues $\lambda_{1}$ and $\lambda_{2}$, then $\Q(\lambda_{1},\lambda_{2})$ is a quadratic field.
	\end{lem}
	\begin{proof}
		Let $A\in \mathrm{GL}_{2}(\Z)$. Then the characteristic polynomial of $A$ is 
		$p(x)=x^{2}-\mathrm{tr}(A)x+\mathrm{det}(A)$ and 
		\begin{equation*}
			\lambda_{\pm}=\frac{\mathrm{tr}(A)\pm \sqrt{(\mathrm{tr}(A))^{2}-4\mathrm{det}(A)}}{2}
		\end{equation*}
		are its eigenvalues. Note that $A$ has equal eigenvalues if and only if $\mathrm{tr}(A)=\pm 2$ and $\mathrm{det}(A)=1$. Since $A$ has distinct eigenvalues and $\mathrm{tr}(A)\neq 0$, we conclude that  the possible values for the trace and the determinant of the matrix $A$ are:
		\begin{enumerate}
			\item[(i)] $\mathrm{tr}(A)=\pm 1$ and $\mathrm{det}(A)=\pm 1$ or
			\item[(ii)] $\mathrm{tr}(A)=\pm 2$ and $\mathrm{det}(A)=-1$ or
			\item[(iii)] $|\mathrm{tr}(A)|>2$ and $\mathrm{det}(A)=\pm 1$.
		\end{enumerate}
		The cases (i) and (ii) are trivial. If $A$ satisfies (iii), then the eigenvalues of $A$ are  irrational. Indeed, suppose that there is an integer $k$ such that $\mathrm{tr}(A)^{2}-4\mathrm{det}(A)=k^{2}$. This gives $\mathrm{tr}(A)^{2}-k^{2}=\pm 4$, which implies $k=\mathrm{tr}(A)\pm 1$ and so $2(\pm \mathrm{tr}(A))-1=\pm 4$, which is a contradiction, because $|\mathrm{tr}(A)|>2$. Therefore, we conclude that $K=\Q(\lambda_{+},\lambda_{-})$ is a quadratic field, i.e.,  there is a square-free integer $d$ such that  $K=\Q(\sqrt{d})$.
	\end{proof}
	
	\begin{lem}\label{Z[lambada]=Ok}
		Let $A\in \mathrm{GL}_{2}(\Z)$ with distinct eigenvalues. If $\mathrm{tr}(A)^{2}-4\mathrm{det}(A)$  is square-free, then $\Z[\lambda]=\mathcal{O}_{K}$, where $K=\Q(\lambda)$ and $\lambda$ is an eigenvalue of $A$.
	\end{lem}
	\begin{proof}
		Note that $\{1,\lambda\}$ is a $\Z$-basis for $\Z[\lambda]$ and that the discriminant of $\{1,\lambda\}$ is $D=\mathrm{tr}(A)^{2}-4\mathrm{det}(A)$. Since $D$ is square-free, it follows that $\{1,\lambda\}$ is a $\Z$-basis for $\mathcal{O}_{K}$ by \cite[Thm. 7.1.8 ]{AW04}. Therefore, $\Z[\lambda]=\mathcal{O}_{K}$.
	\end{proof}

	\begin{prop}[\cite{DF04}, \textsection 16.3, Prop. 14 and Cor. 19]\label{Prop 14 and Cor 19 of DF04}
		The ring of integers $\mathcal{O}_K$ of an algebraic number field $K$ is a Dedekind Domain. In particular, if $I$ is a nonzero ideal of $\mathcal{O}_K$, then every ideal in the quotient $\mathcal{O}_K/I$ is principal. 
	\end{prop}
	
	\begin{defn}
		The \textit{norm} $N(I)$ of an $\mathcal{O}$-ideal $I$ is the index $|\mathcal{O}/I|.$
	\end{defn}
	
	\begin{prop}[\cite{JW09}, Prop. 4.23 and Thm. 4.24] \label{Pro 4.23 and Thm 4.24 of Jw09}
		Let $\mathcal{O}$ be an order in a quadratic field $K$. Then,
		\begin{enumerate}
			\item $I$ is an ideal of $\mathcal{O}$ if and only if $I$ is a $\Z$-module with basis $\{a,b+c\omega\}$, where $\omega$ is as in Proposition \ref{Thm 4.17 of JW09}, $a,b,c\in \Z, \ a>0, \ c>0, \ 0\leq b < a, \ c\mid a, \ c\mid b$, and $ac$ divides the norm of $b+c\omega$.
			\item $N(I)=ac$.
		\end{enumerate}
	\end{prop}
	
	\begin{defn}
		A \textit{fractional ideal} of an order $\mathcal{O}$ is a finitely generated $\mathcal{O}$-submodule $I\neq 0$ of $K$ such that $d I\subseteq \mathcal{O}$ for some $d\in \mathcal{O} \smallsetminus \{0\}$. 
	\end{defn}
	
	The ideals of $\mathcal{O}$ are also fractional ideals of $\mathcal{O}$ (take $d=1$). For clarity these are occasionally called the \textit{integral ideals} of $\mathcal{O}$.
	
	For any $a\in K\smallsetminus \{0\}$ the cyclic $\mathcal{O}$-module $(a):=a\mathcal{O}=\{ax : x\in \mathcal{O}\}$ is called the \textit{principal fractional ideal} generated by $a$.
	
	\begin{defn}\label{equivalent of ideals}
		Two fractional ideals $I$ and $J$ are \textit{equivalent} $(I\sim J)$ if there is $\alpha \in K\smallsetminus \{0\}$ such that $J=\alpha I$.
	\end{defn}
	
	The equivalence established in Definition \ref{equivalent of ideals} is an equivalence relation in the set of all the  fractional ideals of $\mathcal{O}$ and therefore partitions those ideals into distinct ideal classes.
	
	\begin{prop}[\cite{CR62}, Thm. 20.6]\label{Thm 20.6 CR62}
		The number of classes of fractional $\mathcal{O}$-ideals is finite.
	\end{prop}
	
	\begin{defn}
		A fractional $\mathcal{O}$-ideal $I$ is \textit{invertible} if there is another fractional $\mathcal{O}$-ideal $J$ such that $IJ=(\alpha)$, where $\alpha \in \mathcal{O}$. 
	\end{defn}
	
	\begin{prop}[\cite{Neu99}, \textsection 12]
		Let $\mathcal{O}$ be an order. The set $I(\mathcal{O})$ of invertible fractional $\mathcal{O}$-ideals form an abelian group. Moreover, the fractional  principal $\mathcal{O}$-ideals give a subgroup $P(\mathcal{O}) \subseteq I(\mathcal{O})$.
	\end{prop} 
	
	\begin{defn}
		The quotient group $H(\mathcal{O})=I(\mathcal{O})/P(\mathcal{O})$ is called the \textit{ideal class group} (or \textit{Picard group}) of the order $\mathcal{O}$. We denote by $[I]$ the class of an ideal $I \in I(\mathcal{O})$.
	\end{defn}
	
	Since the fractional $\mathcal{O}$-ideals $I$ and $J$ are in the same class of $H(\mathcal{O})$ if and only if $I\sim J$, we have by Proposition \ref{Thm 20.6 CR62} the following corollary.
	
	\begin{cor}
		The ideal class group $H(\mathcal{O})$ is finite.  Its order will be denoted by $h(\mathcal{O})$.
	\end{cor} 
	
	Given an order $\mathcal{O}$  of a quadratic field $K$ with conductor $f$, we say that a nonzero $\mathcal{O}$-ideal $I$ is \textit{prime to} $f$ provided that $I+f\mathcal{O}=\mathcal{O}$. 
	
	\begin{prop}[\cite{Cox89}, Lem. 7.18, Prop. 7.4 and \cite{JW09}, Thm. 4.37]\label{Lem. 7.18 and Prop. 7.4 of Cox89}
		Let $\mathcal{O}$ be an order of conductor $f$.
		\begin{enumerate}
			\item An $\mathcal{O}$-ideal $I$ is prime to $f$ if and only if its norm $N(I)$ is relatively prime to $f$, i.e., $(N(I),f)=1$.
			\item Every $\mathcal{O}$-ideal prime to $f$ is invertible. 
			\item If $I$ is an invertible $\mathcal{O}$-ideal, then there is always some $\mathcal{O}$-ideal $J$ such that $I\sim J$ and $(N(J),f)=1$.
		\end{enumerate}
	\end{prop}
	
	\begin{lem}\label{prime representative for f}
		Let $\mathcal{O}$ be an order of a quadratic field $K$ with conductor $f$. Then each class of ideals in $H(\mathcal{O})$ contains an ideal that is prime to $f$.
	\end{lem}
	\begin{proof}
		Let $[I]\in H(\mathcal{O})$. We can assume that $I$ is an integral ideal of $\mathcal{O}$. Since $I$ is invertible, it follows from Proposition \ref{Lem. 7.18 and Prop. 7.4 of Cox89} (3) that there is an $\mathcal{O}$-ideal $J$ such that $J\sim I$ and $(N(J),f)=1$. By Proposition \ref{Lem. 7.18 and Prop. 7.4 of Cox89}, we have $J$ is an invertible ideal that is prime to $f$ that belongs to the class of $I$.
	\end{proof}
	
	\begin{lem}\label{iso dos aneis estendidos}
		Let $\mathcal{O}$ be an order of a quadratic field $K$ with conductor $f$. For each $\mathcal{O}$-ideal $J$ prime to $f$, the natural ring homomorphism $\mathcal{O}/J\rightarrow \mathcal{O}_{K}/J\mathcal{O}_{K}$ is an isomorphism.
	\end{lem}
	\begin{proof}
		Since $J$ is prime to $f$, we have by definition that $J+f\mathcal{O}=\mathcal{O}$. Thus,
		\begin{equation*}
			J\mathcal{O}_K+f\mathcal{O}_K=(J+f\mathcal{O})\mathcal{O}_K=\mathcal{O}\mathcal{O}_K=\mathcal{O}_K.
		\end{equation*}
		We claim that $J\mathcal{O}_K\cap\mathcal{O}=J$. Indeed, note that
		\begin{eqnarray*}
			J &\subseteq& J\mathcal{O}_K\cap\mathcal{O} \\
			&\subseteq& (J\mathcal{O}_K\cap\mathcal{O})\mathcal{O} \\
			&\subseteq& (J\mathcal{O}_K\cap\mathcal{O})(J+f\mathcal{O}) \\
			&\subseteq& J + f(J\mathcal{O}_K\cap \mathcal{O}) \\
			&\subseteq& J + J \cdot f\mathcal{O}_K \\
			&\subseteq& J + J \ \text{since} \ f\mathcal{O}_K\subseteq \mathcal{O} \\
			&\subseteq& J.
		\end{eqnarray*}
		This gives $J\mathcal{O}_K\cap\mathcal{O}=J$. Now consider natural homomorphism
		\begin{equation*}\label{epimorfismo}
			\varphi:\mathcal{O}\rightarrow \frac{\mathcal{O}_K}{J\mathcal{O}_K}.
		\end{equation*}
		As $J\mathcal{O}_K+f\mathcal{O}_K=\mathcal{O}_K$, we conclude that $\varphi$ is surjective and from $J\mathcal{O}_K\cap\mathcal{O}=J$ we have $\mathrm{ker}(\varphi)= J$. Therefore, the statement follows from the First Isomorphism Theorem for rings.
	\end{proof}
	
	\begin{lem}\label{ideal prin}
		Let $\mathcal{O}$ be an order of a quadratic field. If $I$ is an invertible $\mathcal{O}$-ideal, then $|I/mI|=|\mathcal{O}/m\mathcal{O}|$ for each positive integer $m$.
	\end{lem}
	
	\begin{proof}
		For each $m\in \N$, we have $mI\subseteq I\subseteq \mathcal{O}$. Then
		\begin{equation*}
			\left| \frac{\mathcal{O}}{mI}\right|=\left|\frac{\mathcal{O}}{I} \right|\left|\frac{I}{mI} \right|.
		\end{equation*}
		Thus,
		\begin{equation}\label{eq 1}
			\left|\frac{I}{mI} \right|=\frac{|\mathcal{O}/mI|}{|\mathcal{O}/I|}.
		\end{equation}
		As $mI=m\mathcal{O} I$, we have by \cite[Thm. 4.36]{JW09} that
		\begin{equation}\label{eq 2}
			\left|\frac{\mathcal{O}}{mI} \right|=N(mI)=N(m\mathcal{O})N(I)=\left|\frac{\mathcal{O}}{m\mathcal{O}} \right|\left| \frac{\mathcal{O}}{I}\right|.
		\end{equation}
		Substituting \eqref{eq 2} into \eqref{eq 1} yields
		\begin{equation*}
			\left| \frac{I}{mI}\right|=\left| \frac{\mathcal{O}}{m\mathcal{O}}\right|.
		\end{equation*}
	\end{proof}
	
	Let $I$ be an $\mathcal{O}$-ideal and $m\in \N$. Note that $I/mI$ is an $\mathcal{O}/m\mathcal{O}$-module with the action of $\mathcal{O}/m\mathcal{O}$ on $I/mI$ given by
	\begin{equation*}\label{action11}
		(\theta+m\mathcal{O})\cdot (\alpha+mI):=\theta\alpha +mI.
	\end{equation*}
	This action is well defined because $I$ is an $\mathcal{O}$-module.
	
	\begin{lem}\label{iso local}
		Let $\mathcal{O}$ be an order of a quadratic field $K$ with conductor $f$. If $I$ is an $\mathcal{O}$-ideal prime to $f$, then $$\frac{I}{p^{i}I}\cong\frac{\mathcal{O}}{p^{i}\mathcal{O}}$$ as $\mathcal{O}/p^{i}\mathcal{O}$-modules for each prime number $p$ and $i\in \N$.
	\end{lem}
	\begin{proof}
		We divide the proof in two cases.
		
		\begin{flushleft}
			\textbf{Case 1.} If $p$ does not divide $f$.
		\end{flushleft}
		
		By Proposition \ref{Pro 4.23 and Thm 4.24 of Jw09}, $\{a, b+c\omega\}$ is a $\Z$-basis for $I$ and $N(I)=ac$, where  $a,b, c\in \Z$. Hence, $\{p^{i}a,p^{i}(b+c\omega)\}$ is a $\Z$-basis for $p^{i}I$ and $N(p^{i}I)=p^{2i}ac$. Since $I$ is prime to $f$, we have $(ac, 1)=1$ by Proposition \ref{Lem. 7.18 and Prop. 7.4 of Cox89} (1). This gives $(p^{2i}ac,f)=1$ since $p$ does not divide $f$. Thus,  $p^{i}I$ are invertible ideals prime to $f$ by Proposition \ref{Lem. 7.18 and Prop. 7.4 of Cox89}. Now, we have by Lemma \ref{iso dos aneis estendidos} that 
		\begin{equation*}
			\frac{\mathcal{O}}{p^{i}I}\cong\frac{\mathcal{O}_K}{p^{i}I\mathcal{O}_K}
		\end{equation*}
		as rings. From Proposition \ref{Prop 14 and Cor 19 of DF04}, it follows that $I/p^{i}I$ is a principal ideal of $\mathcal{O}/p^{i}I$. Since 
		\begin{equation*}
			\left|\frac{I}{p^{i}I} \right|=\left|\frac{\mathcal{O}}{p^{i}\mathcal{O}}\right|
		\end{equation*}
		by Lemma \ref{ideal prin} and $I/p^{i}I$ and $\mathcal{O}/p^{i}\mathcal{O}$ are  cyclic $\mathcal{O}/p^{i}\mathcal{O}$-modules, we have
		\begin{equation*}
			\frac{I}{p^{i}I}\cong\frac{\mathcal{O}}{p^{i}\mathcal{O}}
		\end{equation*}
		as $\mathcal{O}/p^{i}\mathcal{O}$-modules.
		
		\begin{flushleft}
			\textbf{Case 2.} If $p$ divides $f$.
		\end{flushleft}
		
		Since $(N(I), f)=(ac,f)=1$ and $p\mid f$,  we have $p\nmid ac$, which implies $p\nmid a$ and $p\nmid c$. Consider the mapping $\varphi: I\rightarrow \mathcal{O}/p^{i}\mathcal{O}$ defined by $\varphi(\alpha)= \alpha a + p^{i}\mathcal{O}$. It is easy to see that $\varphi$ is a group homomorphism. Now let us prove that $\varphi$ is a bijection. 
		\begin{itemize}
			\item $\varphi$ is surjective. 
			
			To prove the statement it is sufficient to prove that $aI+p^{i}\mathcal{O}=\mathcal{O}$. From $(a,p)=1$ we conclude that $(a^2,p^{i})=1$. Note that $\{p^{i},p^{i}\omega\}$ and $\{a^{2}, a(b+c\omega)\}$ are $\Z$-bases for $p^{i}\mathcal{O}$ and $aI$, respectively. As the Diophantine equation  $xa^{2}+yp^{i}=1$ has an integer solution (because $(a^2,p^{i})=1$), we have $1\in aI+p^{i}\mathcal{O}$. Thus, $\mathcal{O}\subseteq aI+p^{i}\mathcal{O}$ and consequently $aI+p^{i}\mathcal{O}=\mathcal{O}$.
			
			\hspace{0.1cm}
			\item $\mathrm{ker}(\varphi)=p^{i}I$.
			
			Remember that  $\{a,b+c\omega\}$ is a $\Z$-basis for $I$ and that $\{1,\omega\}$ is a $\Z$-basis for $\mathcal{O}$. The result is obvious if $a=1$. Suppose $a\neq 1$. Since $p\nmid a$, we have $p^{i}\notin I$.
			
			Let $\beta\in \mathrm{ker}(\varphi)$. Then $\beta a\in p^{i}\mathcal{O}$. So there is $\alpha\in \mathcal{O}$ such that $\beta a= p^{i}\alpha$. We claim that  $\beta\in p^{i}I$. Indeed, we have
			\begin{equation*}
				\beta=ua+v(b+c\omega)  \ \ \text{and} \ \ \alpha=r+s \omega
			\end{equation*}
			where $u,v,r,s\in \Z$. Then,
			\begin{eqnarray*}
				r+s\omega &=&\alpha \\
				&=& \frac{a}{p^{i}}\beta \\
				&=& \frac{a}{p^{i}}(ua+vb) + \frac{a}{p^{i}}vc \omega. \\
			\end{eqnarray*}
			We thus get 
			\begin{equation*}
				r=\frac{a}{p^{i}}(ua+vb) \ \ \text{and} \ \ s=\frac{acv}{p^{i}}.
			\end{equation*}	
			As $p\nmid ac$ and $s\in \Z$, we have $p^{i}\mid v$. Furthermore, from $p\nmid a^{2}$ we conclude that $p^{i}\mid u$, since $r\in \Z$. This gives $u=p^{i}u'$ and $v=p^{i}v'$, where $u',v'\in \Z$. Therefore,
			\begin{equation*}
				\beta=ua+v(b+c\omega)=p^{i}(u'a+v'(b+c\omega)) \in p^{i}I,
			\end{equation*}
			which implies  $\mathrm{ker}(\varphi)\subseteq p^{i}I$. Clearly $p^{i}I \subseteq\mathrm{ker}(\varphi)$ and so $\mathrm{ker}(\varphi)=p^{i}I$.
		\end{itemize}
		
		Therefore, it follows from the previous points that the mapping $\tilde{\varphi}: I/p^{i}I\rightarrow \mathcal{O}/p^{i}\mathcal{O}$ defined by $\tilde{\varphi}(\alpha + p^{i}I)=\alpha a+ p^{i}\mathcal{O}$ is an isomorphism of abelian groups. Since
		
		\begin{eqnarray*}
			\tilde{\varphi}((\theta+p^{i}\mathcal{O})\cdot(\alpha+p^{i}I)) &=& \tilde{\varphi}(\theta\alpha+p^{i}I)  \\
			&=& \theta\alpha a +p^{i}\mathcal{O} \\
			&=& (\theta+p^{i}\mathcal{O})(\alpha a+p^{i}\mathcal{O}) \\
			&=& (\theta+p^{i}\mathcal{O}) \tilde{\varphi}(\alpha+p^{i}I),
		\end{eqnarray*}
		it follows that $\tilde{\varphi}$ is an isomorphism of $\mathcal{O}/p^{i}\mathcal{O}$-modules.
	\end{proof}
	
	\section{Conjugation in $\mathrm{GL}_{2}(\Z)$ and $\mathrm{GL}_{2}(\widehat{\Z})$}\label{Conjugation in GL}
	
	Let $A$ and $B$ be matrices in $\mathrm{GL}_{2}(\Z)$. Consider the semi-direct products $G_{A}=N_{1}\rtimes_{A} C$ and $G_{B}=N_{2}\rtimes_{B} C$, where $N_{1}\cong \Z\times \Z \cong N_{2}$ and $C \cong \Z$.
	
	\begin{lem}\label{iso sse conj}
		Let $A$ and $B$  be matrices in $\mathrm{GL}_{2}(\Z)$ such that none of its eigenvalues is $1$. Then $G_{A}$ and $G_{B}$ are isomorphic if and only if $A$ is conjugate in $\mathrm{GL}_{2}(\Z)$ to $B$ or $B^{-1}$.
	\end{lem}
	\begin{proof}
		It follows from Lemma \ref{condicao2} (i) that if $f:G_{A}\rightarrow G_{B}$ is an isomorphism, then the restriction of $f$ to $N_1$ defines an isomorphism $P: N_1\rightarrow N_2$. Then, $P\in \mathrm{GL}_{2}(\Z)$. Write $f(0,1)=(u,t)$, where $u\in N_2$ and $t\in C$. Note that,
		\begin{equation}\label{Equation 01}
			(0,1)(v,0)(0,-1)=(Av,0).
		\end{equation}
		Computing $f$ in equation \eqref{Equation 01}, we have
		\begin{eqnarray*}
			(PAv,0)&=& f((0,1)(v,0)(0,-1)) \\
			&=& f(0,1)f(v,0)f(0,-1) \\
			&=& (u,t)(Pv,0)(u,t)^{-1} \\
			&=& (u+B^{t}Pv,t)(B^{-t}(-u),-t) \\
			&=& (B^{t}Pv,0).
		\end{eqnarray*}
		Thus, $PAP^{-1}=B^{t}$. Now, let us show that $t=\pm 1$. To this end, note that
		\begin{eqnarray*}
			f(-P^{-1}u,1)&=& f((P^{-1}(-u),0)(0,1)) \\
			&=& (PP^{-1}(-u),0)f(0,1) \\
			&=& (-u,0)(u,t) \\
			&=& (0,t).
		\end{eqnarray*} 
		Since $(0,1)^{t}=(0,t)$, implies that there is $\beta\in G_A$ such that $\beta^{t}=(-P^{-1}u,1)$. Then, $t$ divides $1$ and hence $t=\pm 1$.
		
		Conversely, if there is $P\in \mathrm{GL}_{2}(\Z)$ such that $PAP^{-1}=B$ or $PAP^{-1}=B^{-1}$, then $G_{A}\cong G_{B}$ or $G_{A}\cong G_{B^{-1}}$ by Lemma \ref{prop 2.5}. As $G_{B}\cong G_{B^{-1}}$ by Lemma \ref{G_{A}=G_A^{-1}}, we have $G_{A}\cong G_{B}$ in both cases.
	\end{proof}
	
	The next lemma is a reformulation of Lemma \ref{iso sse conj} for profinite groups.

	\begin{lem}\label{iso sse conj (versao profinita)}
		Let $A$ and $B$  be matrices in $\mathrm{GL}_{2}(\Z)$  such that none of its eigenvalues is $1$. Then $\widehat{G}_{A}$ and $\widehat{G}_{B}$ are isomorphic if and only if $A$ is conjugate in $\mathrm{GL}_{2}(\widehat{\Z})$ to $B$ or $B^{-1}$.
	\end{lem}
	\begin{proof}
		It follows from Lemma \ref{condicao2} (ii) that if $f:\widehat{G}_{A}\rightarrow \widehat{G}_{B}$ is an isomorphism, then the restriction of $f$ to $\widehat{N}_1$ defines an isomorphism $P: \widehat{N}_1\rightarrow \widehat{N}_2$. Then, $P\in \mathrm{GL}_{2}(\widehat{\Z})$. Suppose that $f(0,1)=(u,t)$, where $u\in \widehat{N}_2$ and $t\in \widehat{C}$. Note that,
		\begin{equation*}
			(0,1)(v,0)(0,-1)=(Av,0).
		\end{equation*}
		By similar arguments as in the proof of Lemma \ref{iso sse conj}, we concluded that $PAP^{-1}=B^{t}$ where $t$ is a unit in $\widehat{C}$. Since $\mathrm{det}(B)=\pm 1 \in \Z$ and 
		\begin{equation*}
			\mathrm{det}(B^{t})=(\mathrm{det}(B))^{t} \in \Z,
		\end{equation*}
		we have $t=\pm 1$.	
		
		Conversely, if there is $P\in \mathrm{GL}_{2}(\widehat{\Z})$ such that $PAP^{-1}=B$ or $PAP^{-1}=B^{-1}$, then $\widehat{G}_{A}\cong \widehat{G}_{B}$ or $\widehat{G}_{A}\cong \widehat{G}_{B^{-1}}$ by profinite version of the Lemma \ref{prop 2.5}. We have by Lemma \ref{G_{A}=G_A^{-1}} that $N_2\rtimes_{B^{-1}} C\cong N_{2}\rtimes_{B} C$ and so $\widehat{N}_2\rtimes_{B^{-1}}\widehat{C}\cong \widehat{N}_2\rtimes_{B} \widehat{C}$. Therefore, $\widehat{G}_{A}\cong \widehat{G}_{B}$.
	\end{proof}
	
	\begin{lem}\label{det=det}
		Let $A$ and $B$  be matrices in $\mathrm{GL}_{2}(\Z)$. If $A$ and $B$ are conjugate in $\mathrm{GL}_{2}(\widehat{\Z})$ then,
		\begin{enumerate}
			\item[(i)] $\mathrm{det}(A)=\mathrm{det}(B)$
			\item[(ii)] $\mathrm{tr}(A)=\mathrm{tr}(B)$
			\item[(iii)] $A$ and $B$ have the same characteristic polynomial.
		\end{enumerate}
	\end{lem}
	\begin{proof}
		If $A$ and $B$ are conjugate in $\mathrm{GL}_{2}(\widehat{\Z})$, then for each positive integer $m$, $\mathrm{det}(A)$ is congruent to $\mathrm{det}(B)$ modulo $m$ and $\mathrm{tr}(A)$ is congruent to $\mathrm{tr}(B)$ modulo $m$. Thus, we conclude that $\mathrm{det}(A)=\mathrm{det}(B)$ and $\mathrm{tr}(A)=\mathrm{tr}(B)$.
		
		Since the characteristic polynomials of $A$ and $B$ are $p_{A}(x)=x^{2}-\mathrm{tr}(A)x +\mathrm{det}(A)$ and $p_{B}(x)=x^{2}-\mathrm{tr}(B)x+\mathrm{det}(B)$, (iii) follows from (i) and (ii).
	\end{proof}

	\begin{prop}\label{Obs A and B have the same n of eigenvalues}
		Let $A$ and $B$ be matrices in $\mathrm{GL}_{2}(\Z)$  and consider the semi-direct products $G_A=N_1\rtimes_{A} C$ and $G_B=N_2\rtimes_{B} C$ with $N_1\cong \Z^{2}\cong N_2$ and $C\cong \Z$.  Then, each  isomorphism class in $\mathfrak{g}(G_A)$ contains a group $G_B$ so that $A$ and $B$ have the same eigenvalues.
	\end{prop}
	\begin{proof}
		Let $G_B$ be a representative of any  isomorphism class in $\mathfrak{g}(G_A)$.  Suppose first that $A$ has at least one eigenvalue equal to $1$. Since $\mathrm{det}(A)=\pm 1$, we conclude that the possible eigenvalues of $A$ are or $1$ and $-1$ or all equal to $1$. If the eigenvalues of $A$ are $1$ and $-1$, then $A$ has order $2$ and hence $G_A$ is a $3$-dimensional Bieberbach group with holonomy group of order $2$ (see \cite[Thm. 3.2]{Szc12}). By \cite[Prop. 3.7 and Lem 3.3]{Ner21} we see that $G_B$ is also a $3$-dimensional Bieberbach group with holonomy group of order $2$ and $A$ is conjugate to $B$ in $\mathrm{GL}_{2}(\widehat{\Z})$. Therefore, it follows from Lemma \ref{det=det} that $A$ and $B$ have the same eigenvalues. Now, if $A$ has all the eigenvalues equal to $1$, then $G_{A}$ is a nilpotent group (Lemma \ref{G_A is nipotent iff A has all the eigenvalues equal to 1}). Hence, $\widehat{G}_A$ and $\widehat{G}_B$ are nilpotent groups, which implies that $B$ has all the eigenvalues equal to $1$ (otherwise, the center of $\widehat{G}_{B}$ is trivial, and consequently $\widehat{G}_B$ is not nilpotent).

		Suppose now that none of the eigenvalues of $A$ is $1$. Then none of the eigenvalues of $B$ is $1$, by what has already been proven in the previous paragraph. Hence,  $A$ is conjugate to $B$ or $B^{-1}$  in $\mathrm{GL}_{2}(\widehat{\Z})$ (Lemma \ref{iso sse conj (versao profinita)}). By Lemma \ref{G_{A}=G_A^{-1}} we can assume that $A$ is conjugate to $B$ in $\mathrm{GL}_{2}(\widehat{\Z})$.  Using again Lemma \ref{det=det}, we conclude that $A$ and $B$ have the same eigenvalues.
	\end{proof}
	
	\begin{lem}\label{tr(A^{-1})=-tr(A)}
		Let $A\in \mathrm{GL}_{2}(\Z)$. Then $\mathrm{tr}(A^{-1})=\mathrm{tr}(A)/\mathrm{det}(A)$.
	\end{lem}
	\begin{proof}
		Let $A\in\mathrm{GL}_{2}(\Z)$. It is easy to see that the characteristic polynomial of $A$ has the form $p(x)=x^{2}-\mathrm{tr}(A)x +\mathrm{det}(A)$ and that \begin{equation*}
			\lambda_{\pm}=\frac{\mathrm{tr}(A)\pm \sqrt{(\mathrm{tr}(A))^{2}-4\mathrm{det}(A)}}{2}
		\end{equation*}
		are the eigenvalues. Note that,
		\begin{equation*}
			\frac{\mathrm{det}(A)}{\lambda_{+}}=\lambda_{-}.
		\end{equation*}
		Since the trace of $A^{-1}$ is equal to the trace of the diagonal matrix $A^{-1}$, we have
		\begin{eqnarray*}
			\mathrm{tr}(A^{-1}) &=& \frac{1}{\lambda_{+}}+\frac{1}{\lambda_{-}} \\
			&=& \frac{1}{\lambda_{+}}+ \frac{\lambda_{+}}{\mathrm{det}(A)} \\
			&=& \frac{2}{\mathrm{tr}(A)+\sqrt{(\mathrm{tr}(A))^{2}-4\mathrm{det}(A)}} + \frac{\mathrm{tr}(A)+\sqrt{(\mathrm{tr}(A))^{2}-4\mathrm{det}(A)}}{2\mathrm{det}(A)} \\
			&=& \frac{\mathrm{tr}(A)-\sqrt{(\mathrm{tr}(A))^{2}-4\mathrm{det}(A)}}{2\mathrm{det}(A)} + \frac{\mathrm{tr}(A)+\sqrt{(\mathrm{tr}(A))^{2}-4\mathrm{det}(A)}}{2\mathrm{det}(A)}  \\
			&=& \frac{\mathrm{tr}(A)}{\mathrm{det}(A)}.
		\end{eqnarray*}
	\end{proof}
	
	\begin{cor}
		Let $A\in\mathrm{GL}_{2}(\Z)$ such that $\mathrm{det}(A)=-1$ and none of its eigenvalues is $1$. Then the cardinality of the genus $\mathfrak{g}(G_{A})$ is exactly the number of conjugacy classes of matrices in $\mathrm{GL}_{2}(\Z)$  in the conjugacy class of $A$ in $\mathrm{GL}_{2}(\widehat{\Z})$.
	\end{cor}
	\begin{proof}
		Let $G_B\in \mathfrak{g}(G_A)$. By Proposition \ref{Obs A and B have the same n of eigenvalues}  none of the eigenvalues of $B$ is $1$. From Lemma \ref{iso sse conj (versao profinita)}, implies that $A$ is conjugate in $\mathrm{GL}_{2}(\widehat{\Z})$ to $B$ or $B^{-1}$. By Lemma \ref{G_{A}=G_A^{-1}} we can assume that $A$ is conjugate to $B$ in $\mathrm{GL}_{2}(\widehat{\Z})$. Thus $\mathrm{det}(A)=\mathrm{det}(B)$ and $\mathrm{tr}(A)=\mathrm{tr}(B)$ by Lemma \ref{det=det}. Since $\mathrm{det}(A)=-1$, we have $\mathrm{tr}(B^{-1})=-\mathrm{tr}(B)$ (Lemma \ref{tr(A^{-1})=-tr(A)}). Therefore, the statement follows by Lemmas \ref{iso sse conj} and \ref{iso sse conj (versao profinita)}.
	\end{proof}
	
	Thus, if $ G_B\in \mathfrak{g}(G_A)$, we can assume that $A$ and $B$ have the same characteristic polynomial. Therefore, to study the extent to which matrices $A$ and $B$ are conjugate in $\mathrm{GL}_{2}(\Z)$, we need to study the set of conjugacy classes in $\mathrm{GL}_{2}(\Z)$ of matrices which have the same characteristic polynomial. For this purpose, the next theorem tells us that it is enough to study the ideal classes of orders of quadratic fields.
	
	\begin{prop}[Latimer-MacDuffee, \cite{New72}, Thm. III. 13]\label{Latimer-MacDffe}
		Let $f(x)\in \Z[x]$ be a monic polynomial of degree $n$ that is irreducible over $\Q$ and let $\lambda$ be a root of $f(x)$. Then there is a one-to-one correspondence between the $\Z$-similarity classes of $n \times n$ matrices over $\Z$ with characteristic polynomial $f$ and the ideal classes of the order $\Z[\lambda]$.
	\end{prop}
	
	Here we are interested in the case in which $n=2$. Let
	\begin{equation*}
		A=\begin{pmatrix}
			a & b \\ c  & d
		\end{pmatrix}\in \mathrm{GL}_{2}(\Z)
	\end{equation*}
	and let $\lambda$ be an eigenvalue of $A$. Consider the $\Z$-module 
	\begin{equation}\label{I_A}
		I_A=\langle b, \lambda-a \rangle = \{mb+n(\lambda-a):m,n\in \Z\}.
	\end{equation}
	
	\begin{prop}[\cite{Hen97}, Thm. 1]
		Let $A$ be a matrix in $\mathrm{ GL}_{2}(\Z)$. Then,	$I_A$ is an ideal of the order $\Z[\lambda]$.
	\end{prop}
	
	We have the following version of Proposition \ref{Latimer-MacDffe} for $\mathrm{GL}_{2}(\Z)$.
	
	\begin{prop}[\cite{Hen97}, Thm. 2]
		Let $A$ and $B$ be matrices in $\mathrm{GL}_{2}(\Z)$ with the same characteristic polynomial $f(x)$. Then, $A$ and $B$ are conjugate matrices in $\mathrm{GL}_{2}(\Z)$ if and only if the corresponding ideals $I_A$ and $I_B$ are in the same ideal class of $ \Z[\lambda]$, where $\lambda$ is a root of $f(x)$.
	\end{prop}
	
	Let $\lambda$ be an eigenvalue of the matrix $A\in \mathrm{GL}_{2}(\Z)$. For any fractional $\Z[\lambda]$-ideal  $I$, the multiplication by $\lambda$ is a $\Z$-linear map, i.e., $m_{I,\lambda}: I\rightarrow I $ defined by $x\mapsto x \lambda$. Since $I$ is a $\Z$-module, choosing a $\Z$-basis we can represent $m_{I,\lambda}$ by a $2\times 2$ matrix  $[m_{I,\lambda }]$ over $\Z$. Note that, if $\lambda$ is a unit of $\Z[\lambda]$, then $m_{I,\lambda}$ is a bijection (because the mapping $m_{I,\lambda^{- 1}}:I\rightarrow I$ given by $x\mapsto x\lambda^{-1}$ is an inverse to $m_{I,\lambda}$). Thus, if $\lambda$ is a unit of $\Z[\lambda]$, then $[m_{I,\lambda}]\in \mathrm{GL}_{2}(\Z)$.
	
	\begin{lem}\label{matrizes nao conjugadas}
		Let $A\in\mathrm{GL}_{2}(\Z)$ with distinct eigenvalues and $\mathrm{tr}(A)\neq 0$. Let $I$ and $J$ be ideals in distinct classes  of $H(\Z[\lambda])$, where $\lambda$ is an eigenvalue of $A$. Then, the matrices $[m_{I,\lambda}]$ and $[m_{J,\lambda}]$ are not conjugate in $\mathrm{GL}_{2}(\Z)$.
	\end{lem}
	\begin{proof}
		Initially, since $\lambda$ is a root of the polynomial $p(x)=x^{2}-\mathrm{tr}(A)x+\mathrm{det}(A)$, we have that $\lambda$ is an algebraic integer of the quadratic field $K=\Q(\lambda)$. As $\mathrm{det}(A)=\pm 1$,  implies that  the norm of $\lambda$ is equal to $\pm 1$ and hence  $\lambda$ is a unit of the order $\Z[\lambda]$.
		
		Let $B_1=\{v_1,w_1\}$ and $B_2=\{v_2,w_2\}$ be fixed $\Z$-bases  for the $\Z$-modules $I$ and $J$, respectively. Consider the matrix representations $[m_{I,\lambda}]$ and $[m_{J,\lambda}]$ of the $\Z$-linear automorphisms $m_{I,\lambda}:I\rightarrow I$ and $m_{J,\lambda}:J\rightarrow J$ defined by the multiplications by $\lambda$ with respect to the bases $B_1$ and $B_2$, respectively. Let us denote by $\{e_1,e_2\}$ the canonical basis of $\Z^{2}$. Now, consider the $\Z$-linear isomorphisms $f:I\rightarrow \Z^{2}$ and $g:J\rightarrow \Z^{2}$ such that $f(v_1)=e_1=g( v_2)$ and $f(w_1)=e_2=g(w_2)$. It follows by construction that the diagrams
		\[
		\begin{CD}
			I @>{f}>> \Z^{2} \\
			@VV{m_{I,\lambda}}V @VV{[m_{I,\lambda}]}V \\
			I @>{f}>> \Z^{2}
		\end{CD} \hspace{1cm}\text{and} \hspace{1cm}
		\begin{CD}
			J @>{g}>> \Z^{2} \\
			@VV{m_{J,\lambda}}V @VV{[m_{J,\lambda}]}V \\
			J @>{g}>> \Z^{2}
		\end{CD}
		\]
		are commutative.
		
		Suppose that the matrices $[m_{I,\lambda}]$ and $[m_{J,\lambda}]$ are conjugate in $\mathrm{GL}_{2}(\Z)$, i.e., there is  $P\in\mathrm{GL}_{2}(\Z)$ such that $P^{-1}[m_{I,\lambda}]P=[m_{J,\lambda}]$. Thus, we obtain the following diagram of $\Z$-linear isomorphisms
		\[
		\begin{CD}
			I @>{f}>> \Z^{2} @>{P}>> \Z^{2} @>{g^{-1}}>> J  \\
			@VV{m_{I,\lambda}}V @VV{[m_{I,\lambda}]}V @VV{[m_{J,\lambda}]}V @VV{m_{J,\lambda}}V\\
			I @>{f}>> \Z^{2} @>{P}>> \Z^{2} @>{g^{-1}}>> J \ .
		\end{CD}
		\]
		As each square in the diagram is commutative, we have that the whole diagram is commutative. We claim that the composite map $g^{-1}Pf:I\rightarrow J$ is an $\Z[\lambda]$-module isomorphism. Indeed, note that for every $\alpha\in I$ we have
		\begin{equation*}
			g^{-1}Pf(\alpha\lambda)=g^{-1}Pf(m_{I,\lambda}(\alpha))=m_{J,\lambda}g^{-1}Pf(\alpha)=g^{-1}Pf(\alpha)\lambda.
		\end{equation*}
		Since the composite map $g^{-1}Pf:I\rightarrow J$ is clearly a $\Z$-linear isomorphism, it follows that $g^{-1}Pf:I\rightarrow J$ is an $\Z[\lambda]$-module isomorphism. Now consider an extension of $g^{-1}Pf:I\rightarrow J$ to a $K$-isomorphism $\psi: KI\rightarrow KJ$. Thus, we have
		\begin{equation*}
			\psi(m\xi)=\psi(m)\xi, \ \  \ m\in KI, \ \xi\in K.
		\end{equation*}
		In particular, if $m=1$ and restricting $\xi$ to $I$, we have 
		\begin{equation*}
			g^{-1}Pf(\xi)=\psi(\xi)=\psi(1)\xi, \ \ \xi\in I.
		\end{equation*}
		Therefore, $g^{-1}Pf:I\rightarrow J$ is the multiplication by $\psi(1)$ and so $I=\psi(1)J$. Then $I$ and $J$ are in the same class of ideals, a contradiction. Therefore, the matrices $[m_{I,\lambda}]$ and $[m_{J,\lambda}]$ are not conjugate in $\mathrm{GL}_{2}(\Z)$.
	\end{proof}
	
	Let $A\in \mathrm{GL}_{2}(\Z)$ be a fixed matrix. We define
	\begin{equation*}
		\approx_A\text{-class}:=\{B\in\mathrm{GL}_{2}(\Z): B \ \text{is conjugate to}\ A \ \text{in} \ \mathrm{GL}_{2}(\widehat{\Z})\}.
	\end{equation*}
	
	\begin{lem}\label{existence of at least h(K) not coj}
		Let $A$ be a matrix in $\mathrm{GL}_{2}(\Z)$ such that its eigenvalues are distinct and $\mathrm{tr}(A)\neq 0$. Let $\lambda$ be an eigenvalue of $A$. If the conjugacy class of $A$ corresponds to a class $[I_A]$ of an invertible ideal $I_A$ of $\Z[\lambda]$, where $I_A$ is as in \eqref{I_A}, then there are at least $h (\Z[\lambda])$ matrices in $\approx_A$-class that are not conjugate to each other.
	\end{lem}
	\begin{proof}
		Let $I_1,\cdots, I_k$ be ideals in  distinct classes of the group $H(\mathcal{\Z[\lambda]})$. By Lemma \ref{prime representative for f}, we can assume that each $I_i$ is  prime to the conductor $f=|\mathcal{O}_K:\Z[\lambda]|$. Let $B_i$ be a fixed $\Z$-basis for each $I_i$ and consider $A_i\in \mathrm{GL}_{2}(\Z)$ the matrix of multiplication by $\lambda$ with respect to the basis $B_i$, $i=1,\cdots,k$. By Lemma \ref{matrizes nao conjugadas}, $A_i$ is not conjugate to $A_j$ in $\mathrm{GL}_{2}(\Z)$, whenever $i\neq j$.
		
		On the other hand, by Lemma \ref{iso local}
		\begin{equation}\label{equation1}
			\frac{I_i}{p^{l}I_{i}}\cong\frac{\Z[\lambda]}{p^{l}\Z[\lambda]}
		\end{equation}
		as $\Z[\lambda]/p^{l}\Z[\lambda]$-modules, for every prime $p$ and $l\in \N$. Now consider the canonical images $\bar{A}_i,\bar{A}_1$ in $\mathrm{GL}_{2}(\Z/p^{l}\Z)$ of the matrices $A_i, A_1$. From \eqref{equation1}, we see that $\bar{A}_i$ is conjugate to $\bar{A}_1$ in $\mathrm{GL}_{2}(\Z/p^{l}\Z)$ for every $l\in \N$ and  $i=1,\cdots, k$. Hence, the matrices $A_i$ and $A_1$ are conjugate in $\mathrm{GL}_{2}(\mathbf{Z}_p)$ for every prime $p$ and $i=1,\cdots,k$. Consequently, $A_i$ is conjugate to $A_1$ in $$\mathrm{GL}_{2}(\widehat{\Z})=\prod_{p}\mathrm{GL}_{2}(\mathbf{Z}_p),$$ for each $i=1,\cdots,k$.
		
		Since the conjugacy class of the matrix $A$ corresponds to a class $[I_A]$ of an invertible ideal $I_A$ of $\Z[\lambda]$, it follows by Proposition \ref{Latimer-MacDffe} that there is $j\in\{ 1,\cdots,k\}$ such that $A$ is conjugate to $A_j$ in $\mathrm{GL}_{2}(\Z)$. Therefore, there are at least $h(\Z[\lambda])$ matrices in $\approx_A$-class that are not conjugate to each other.
	\end{proof}

	The next result will be useful to us, later on.

	\begin{lem}[\cite{CT91}, Lem 3.2]\label{Lem 3.2}
		Let $A$ and $B$  be matrices in $\mathrm{GL}_{2}(\Z)$. If $\mathrm{det}(A)=\mathrm{det}(B)$ and $\mathrm{tr}(A)=\mathrm{tr}(B)\neq 0$, then $A$ is conjugate to $B$ in $\mathrm{ GL}_{2}(\Z)$ if and only if $A^{2}$ is conjugate to $B^{2}$ in $\mathrm{ GL}_{2}(\Z)$.
	\end{lem}

	\section{Proof of main results}\label{main results}

	\begin{proof}[Proof of Theorem \ref{Caso nilpotente Th}]
		If $A=I$, then the group $G_A\cong \Z^{3}$ and the result follows from \cite[Prop. 3.1]{Reid18}. Now suppose  $A\neq I$ and that its eigenvalues are all equal to $1$. Let $B\in \mathrm{GL}_{2}(\Z)$ such that $\widehat{G}_A\cong \widehat{G}_B$. By Proposition \ref{Obs A and B have the same n of eigenvalues} we can assume that $B$ also has all its eigenvalues equal to $1$. By Lemma \ref{G_A is nipotent iff A has all the eigenvalues equal to 1}, $G_A$ and $G_B$ are nilpotent groups with nilpotency class $2$. Also, note that $G_A$ and $G_B$ has  Hirsch number $3$. Since $G_A$ and $G_B$ are finitely generated and torsion-free, it follows by Proposition \ref{Grunewald-Scharlau} that $G_A\cong G_B$, and  $\#\mathfrak{g}(\mathcal{PF},G_A)=1$ as claimed.
	\end{proof}
	
	\begin{proof}[Proof of Theorem \ref{Caso nao nilpotente Th}]
		Let $A$ and $B$ be matrices in $\mathrm{GL}_{2}(\Z)$ and let $\lambda$ be an eigenvalue of $A$. Consider the semi-direct products $G_{A}=N_{1}\rtimes_{A}\Z$ and $G_{B}=N_{2}\rtimes_{B}\Z$, where $N_{1}\cong \Z\times\Z\cong N_{2}$ and suppose that
		\begin{eqnarray*}
			\widehat{G}_{B}\cong\widehat{G}_{A}.
		\end{eqnarray*}
		
		\begin{enumerate}
			\item[(i)] Let us assume that all eigenvalues of $A$ are distinct, $\mathrm{tr}(A)\neq 0$ and that the conjugacy class of $A$ in $\mathrm{GL}_{2}(\Z)$ corresponds to a class $[I_A]$ of invertible ideals of $\Z[\lambda]$. Then, none of the eigenvalues of $A$ are equal to $1$, and since $\widehat{G}_A\cong \widehat{G}_B$, we can assume that $A$ and $B$ have the same eigenvalues by Proposition \ref{Obs A and B have the same n of eigenvalues}. By Lemma \ref{existence of at least h(K) not coj} there are at least $h(\Z[\lambda])$ matrices in $\approx_A$-class that are not conjugate to each other in $\mathrm{ GL}_{2}(\Z)$ but  are conjugate in $\mathrm{GL}_{2}(\widehat{\Z})$. We now  consider two cases to establish the lower bound.
			
			\begin{flushleft}
				\textbf{Case (a)} If $\mathrm{det}(A)=-1$.
			\end{flushleft}
			
			Note that, if $A,B\in \approx_A$-class we have  $\mathrm{det}(A)=\mathrm{det}(B)=-1$ and $\mathrm{tr}(A)= \mathrm{tr}(B)$ by Lemma \ref{det=det}. Thus, $\mathrm{tr}(B^{-1})=-\mathrm{tr}(B)$ by Lemma \ref{tr(A^{-1})=-tr(A)} and hence $B^{-1 }$ is not conjugate to any matrix of $\approx_A$-class. Therefore, $$h(\Z[\lambda])\leq \#\mathfrak{g}(\mathcal{PF}, G_A)$$ by Lemmas \ref{iso sse conj} and \ref{iso sse conj (versao profinita)}.
			
			\begin{flushleft}
				\textbf{Case (b)} If $\mathrm{det}(A)=1$.
			\end{flushleft}
			
			If there are one or two conjugacy classes in $\approx_A$-class, then the result immediately follows from Lemma \ref{iso sse conj}. Now suppose that $\approx_A$-class has more than two matrices that are not conjugate to each other and that $A,B,C\in \mathrm{GL}_{2}(\Z)$ are three of these matrices. Note that, if $A^{-1}$ is conjugate to $B$ in $\mathrm{GL}_{2}(\Z)$ and $C^{-1}$ is also conjugate to $B$ in $\mathrm{GL}_{2}(\Z)$, then $A^{-1}$ is conjugate to $C^{-1}$ in $\mathrm{GL}_{2}(\Z)$, which implies that $A$ is conjugate to $C$ in $\mathrm{GL}_{2}(\Z)$,  a contradiction. Therefore, $A^{-1}$ and $C^{-1}$ cannot be conjugate to more than one matrix of the subset of $\approx_A$-class formed by the matrices that are not conjugate to each other. Thus, by Lemmas \ref{iso sse conj} and \ref{iso sse conj (versao profinita)} we have
			\begin{equation*}
				h(\Z[\lambda])/2\leq \#\mathfrak{g}(\mathcal{PF},G_A).
			\end{equation*}
			
			On the other hand, by Proposition \ref{Obs A and B have the same n of eigenvalues} we can assume that $A$ and $B$ have the same characteristic polynomial $p_{A}(x)$. Now, Theorem \ref{Latimer-MacDffe} tells us that there is a one-to-one correspondence between the  conjugacy classes of matrices in $\mathrm{GL}_{2}(\Z)$ that have characteristic polynomial $p_A(x)$ and the ideal classes of the order $\Z[\lambda]$. Hence, if $\tilde{h}(\lambda)$ denotes the number of ideal classes of the order $\Z[\lambda]$, it follows by Lemma \ref{iso sse conj} that $$\# \mathfrak{g}(\mathcal{PF},G_A)\leq \tilde{h}(\lambda).$$ Therefore,
			\begin{equation*}
				\left\{
				\begin{array}{c c}
					h(\Z[\lambda])\leq \#\mathfrak{g}(\mathcal{PF},G_A)\leq \tilde{h}(\lambda), & \text{if} \ \mathrm{det}(A)=-1\\
					h(\Z[\lambda])/2\leq \#\mathfrak{g}(\mathcal{PF},G_A)\leq \tilde{h}(\lambda), & \text{if} \ \mathrm{det}(A)=1.	
				\end{array}
				\right.
			\end{equation*}
			
			\item[(ii)] We divide the proof into two cases.
			
			\begin{flushleft}
				\textbf{Case 1.} If $A$ has all its eigenvalues equal to $-1$.
			\end{flushleft}
			
			It is easy to see that the groups $G_{B^{2}}=N_{2}\rtimes_{B^{2}}\Z$ and $G_{A^{2}}=N_{1}\rtimes_{A^{2}}\Z$ are subgroups of $G_{B}$ and $G_{A}$ (corresponding to the replacement of the second factor $\Z$ by $2\Z$). Thus,  the following sequences $$1\rightarrow G_{B^{2}}\rightarrow G_{B}\rightarrow C_{B}\rightarrow 1 \ \ \text{and} \ \ 1\rightarrow G_{A^{2}}\rightarrow G_{A}\rightarrow C_{A}\rightarrow 1$$ are exact, where $C_{B}\cong \Z/2\Z\cong C_{A}$. Note that the epimorphism $p_{i}:G_{i}\rightarrow C_{i}$ induces an epimorphism of profinite groups $\widehat{p}_{i}:\widehat{G}_{i}\rightarrow C_{i}$ where  $\mathrm{ker} (\widehat{p}_{i})=\widehat{G}_{i^{2}}$, for $i \in \{A,B\}$. Hence the sequences 
			\begin{equation*}
				1\rightarrow \widehat{G}_{B^{2}}\rightarrow \widehat{G}_{B}\stackrel{\widehat{p}_{B}}{\rightarrow} C_{B}\rightarrow 1 \ \ \text{and} \ \ 1\rightarrow \widehat{G}_{A^{2}}\rightarrow \widehat{G}_{A}\stackrel{\widehat{p}_{A}}{\rightarrow} C_{A}\rightarrow 1
			\end{equation*}
			are exact.
			
			Since none of the eigenvalues of the matrices $A$ and $B$ are equal to $1$  (Proposition \ref{Obs A and B have the same n of eigenvalues}), we have by Lemma \ref{condicao2} that the restriction of any isomorphism $f:\widehat{G}_{B}\rightarrow \widehat{G}_{A}$ to $\widehat{N}_{2}$ is an isomorphism onto $\widehat{N}_{1}$. Hence $f$ induces an isomorphism $\sigma: C_{B}\rightarrow C_{A}$. Thus, we get the following commutative diagram where the restriction of map $f$ for $\widehat{G}_{B^{2}}$ is an isomorphism from  $\widehat{G}_{B^{2}}$ to  $\widehat{G}_{A^{2}}$
			\[
			\begin{CD}
				1 @>{}>> \widehat{G}_{B^{2}} @>{}>> \widehat{G}_{B} @>{\widehat{p}_{B}}>> C_{B} @>>> 1 \\
				@. @VV{}V @VV{f}V @VV{\sigma}V\\
				1 @>{}>> \widehat{G}_{A^{2}} @>{}>> \widehat{G}_{A} @>{\widehat{p}_{A}}>> C_{A} @>>> 1 \ .
			\end{CD}
			\] 
			Note that $A^{2}$ has all its eigenvalues equal to 1, hence there is an isomorphism $\phi: G_{B^{2}}\rightarrow G_{A^{2}}$  by Theorem \ref{Caso nilpotente Th}. Since $N_{2}\rtimes_{B} 2\Z = G_{B^{2}}$ and $N_{1}\rtimes_{A} 2\Z = G_{A^{2}}$ and none of the eigenvalues of matrices $B$ and $A$ are equal to $1$, it follows from Lemma \ref{condicao2} that $\phi(N_{2})=N_{1}$. We conclude from Lemma \ref{cor 2.2} that the matrices $B^{2}$ and $A^{2}$ are conjugate in $\mathrm{GL}_{2}(\Z)$. 
			Again, we can assume that $A$ and $B$ are conjugate in $\mathrm{GL}_{2}(\widehat{\Z})$ (see Lemmas \ref{G_{A}=G_A^{-1}} and \ref{iso sse conj (versao profinita)}). We thus get $\mathrm{tr}(B)=\mathrm{tr}(A)$ and  $\mathrm{det}(B)=\mathrm{det}(A)$ (Lemma \ref{det=det}). Lemma \ref{Lem 3.2} now shows that $B$ and $A$ are conjugate in $\mathrm{GL}_{2}(\Z)$. Hence $G_{B}\cong G_{A}$ (Lemma \ref{iso sse conj}). Therefore,  $\#\mathfrak{g}(\mathcal{PF},G_{A})=h(\Z)=1$, since $\Z$ is a principal ideal domain. 
			
			\begin{flushleft}
				\textbf{Case 2.} If  $\mathrm{tr}(A)=0$.
			\end{flushleft}
			
			First, we assume that $\mathrm{det}(A)=1$. In this case, the eigenvalues of $A$ are distinct and different from $1$. From the proof of Proposition \ref{Obs A and B have the same n of eigenvalues}, we see that none of the eigenvalues of $B$ is $1$. As $\widehat{G}_A\cong\widehat{G}_B$, it follows by  Lemmas \ref{G_{A}=G_A^{-1}} and \ref{iso sse conj (versao profinita)} that we can assume that $A$ and $B$ are conjugate in $\mathrm{ GL}_{2}(\widehat{\Z})$. This gives $\mathrm{tr}(B)=\mathrm{tr}(A)=0$ and $\mathrm{det}(B)=\mathrm{det}(A)=1$ (Lemma \ref{det=det}). Now, by \cite[Thm. 3]{Hen97} we see that every matrix with  trace equal to zero and determinant equal to $1$ is conjugate in $\mathrm{GL}_{2}(\Z)$ to $$\begin{pmatrix}
				0 & -1 \\ 1  & 0
			\end{pmatrix}.$$Hence, the matrices $A$ and $B$ are conjugate in $\mathrm{GL}_{2}(\Z)$ and so $G_{A}\cong G_B$ by the Lemma \ref{iso sse conj}. Therefore, $\#\mathfrak{g}(\mathcal{PF},G_A)=1$. On the other hand, since the characteristic polynomial of $A$ is $p(x)=x^2+1$, we have that $\sqrt{-1}$ and $-\sqrt{-1}$ are the eigenvalues of $A$. From \cite[Thm. 5.4.2]{AW04} we see that $\Z[\sqrt{-1}]$ is the ring of integers of the field $\Q(\sqrt{-1})$. Thus, $\#\mathfrak{g}(\mathcal{PF},G_A)=h(\Z[\sqrt{-1}])=1$ by \cite[p. 325]{AW04}.
			
			Finally, suppose that $\mathrm{det}(A)=-1$. By \cite[Thm. 3]{Hen97} we have that $A$ is conjugate in $\mathrm{GL}_{2}(\Z)$ to $$P=\begin{pmatrix}
				0 & 1 \\ 1  & 0
			\end{pmatrix}.$$ Since $P^2=I$, it follows that $G_A$ is a $3$-dimensional Bieberbach group with holonomy group $\Z/2\Z$ (see \cite[Thm. 3.2]{Szc12}). Therefore, by \cite[Cor. 1.2]{Ner21}  we have  $\#\mathfrak{g}(\mathcal{PF},G_A)=1$. In this case, the characteristic polynomial of $A$ is $p(x)=x^{2}-1$, so $1$ and $-1$ are the eigenvalues of $A$. Thus, as in Case 1, we see that $\#\mathfrak{g}(\mathcal{PF},G_A)=h(\Z)=1$.
		\end{enumerate}
	\end{proof}

	\begin{proof}[Proof of Corollary \ref{cor1}]
		By Lemma \ref{Z[lambada]=Ok}, we have $h(\Z[\lambda])=\tilde{h}(\lambda)$. Therefore, $\#\mathfrak{g}(\mathcal{PF},G_A)=h(\Z[\lambda])$ by Theorem \ref{Caso nao nilpotente Th}.
	\end{proof}
	
	\begin{proof}[Proof of Corollary \ref{cor2}]
		This is an immediate consequence of the Theorems \ref{Caso nilpotente Th} and \ref{Caso nao nilpotente Th}.
	\end{proof}
	
	\section{Examples}\label{examples}
	
	In this section we give some examples of applications of the Theorems \ref{Caso nilpotente Th} and \ref{Caso nao nilpotente Th}.
	
	\begin{example}
		Consider \begin{equation*}
			A=\begin{pmatrix}
				2 & 1 \\ 5  & 2
			\end{pmatrix} \in \mathrm{GL}_{2}(\Z)
		\end{equation*}
		with characteristic polynomial $f_A(x)=x^2-4x-1$. So the eigenvalues of $A$ are $\lambda_{\pm}=2\pm \sqrt{5}.$
		Since $5\equiv 1 \ (\mathrm{mod} \ 4)$, it follows by \cite[Thm. 5.4.2]{AW04} that the order $\Z[\lambda_{+}]=\Z[\sqrt{ 5}]$ is not maximal, i.e., $\Z[\sqrt{5}]$ is not the ring of integers of the quadratic field $\Q(\sqrt{5})$. Consider $I=\langle 2, 1+\sqrt{5}\rangle$. By Proposition \ref{Pro 4.23 and Thm 4.24 of Jw09} (1),  we see that $I$ is an ideal of $\Z[\sqrt{5}]$. By \cite[Ex. 1.32]{Mol01},  $I$ is not invertible in $\Z[\sqrt{5}]$. This implies that the order $h(\Z[\sqrt{5}])$ of the ideal class group of $\Z[\sqrt{5}]$ is strictly less than the number of ideal classes $\tilde{h}(\lambda_{+})$ of $\Z[\sqrt{5}]$. As the corresponding ideal to the matrix $A$ is $I_A=\langle 1, \lambda_{+}-2\rangle= \Z[\sqrt{5}]$ (see \eqref{I_A}), which is invertible, and $\mathrm{det}(A)=-1$, it follows by Theorem \ref{Caso nao nilpotente Th} that holds the following inequality
		\begin{equation*}
			h(\Z[\sqrt{5}])\leq \#\mathfrak{g}(\mathcal{PF},\Z^{2}\rtimes_{A} \Z)\leq \tilde{h}(\lambda_{+}) \ \  \text{with} \ \ h(\Z[\sqrt{5}])< \tilde{h}(\lambda_{+}).
		\end{equation*}
	\end{example}
	
	The next example shows that torus bundles modeled on the geometry $\mathbf{Sol}$ are not determined by the profinite completions of their fundamental groups.
	
	\begin{example}
		Consider \begin{equation*}
			B=\begin{pmatrix}
				4 & 7 \\ 3  & 5
			\end{pmatrix} \in \mathrm{GL}_{2}(\Z).
		\end{equation*} 
		Note that the characteristic polynomial of $B$ is $f_B(x)=x^2-9x-1$ and that
		\begin{equation*}
			\lambda_{\pm}=\frac{9\pm \sqrt{85}}{2}
		\end{equation*}
		are the eigenvalues. Since $85$ is  square-free and $\mathrm{det}(B)=-1$, it follows from Corollary \ref{cor1} that
		\begin{equation*}
			\#\mathfrak{g}(\mathcal{PF},\Z^{2}\rtimes_{B}\Z)=h(\Z[\lambda_{+}]).
		\end{equation*}
		Now, by \cite[Table 8]{AW04} we see that $h(\Z[\lambda_{+}])=2$. As $|\mathrm{tr}(B)|>2$, it follows that the torus bundle $M_B$, such that $\pi_1(M_B)\cong \Z^{2}\rtimes_{B} \Z $, is modeled on the geometry $\mathbf{Sol}$ (Cf. \cite[Thm. 5.5]{Sco83}).
	\end{example}
	
	In contrast, Corollary \ref{cor2} tells us that there is a family of torus bundles, modeled on the geometry $\mathbf{Sol}$, that are determined by the profinite completions of their fundamental groups.
	
	\begin{example}
		Consider the matrices
		\begin{equation*}
			C=\begin{pmatrix}
				2 & 7 \\ 1  & 3
			\end{pmatrix}, \ D=\begin{pmatrix}
				5 & 1 \\ 11  & 2
			\end{pmatrix}\in \mathrm{GL}_{2}(\Z).
		\end{equation*}	
		Since the trace of these matrices are absolutely greater than $2$, we have that the torus bundles $M_C$ and $M_D$, such that
		\begin{equation*}
			\pi_1(M_C)\cong \Z^{2}\rtimes_{C}\Z \ \text{and} \ \pi_1(M_D)\cong \Z^{2}\rtimes_{D}\Z,
		\end{equation*}
		are all modeled on the geometry $\mathbf{Sol}$ (Cf. \cite[Thm. 5.5]{Sco83}). Now, note that the characteristic polynomials of $C$ and $D$ are
		\begin{equation*}
			f_{C}(x)=x^2-5x-1, \ f_D(x)=x^2-7x-1,
		\end{equation*}
		respectively. Hence, the eigenvalues of $C$ and $D$ are
		\begin{equation*}
			\lambda_{\pm,C}=\frac{5\pm\sqrt{29}}{2}, \ \lambda_{\pm,D}=\frac{7\pm\sqrt{53}}{2},
		\end{equation*}
		respectively. Since $29$ and $53$ are square-free numbers and $\mathrm{det}(C)=\mathrm{det}(D)=-1$, it follows from Corollary \ref{cor1} with \cite[Table 8]{AW04} that
		\begin{equation*}
			\#\mathfrak{g}(\mathcal{PF},\Z^{2}\rtimes_{C}\Z)=1=  \#\mathfrak{g}(\mathcal{PF},\Z^{2}\rtimes_{D}\Z).
		\end{equation*}
		Therefore, the torus bundles $M_C$ and $M_D$ are determined among  $3$-manifolds by the profinite completion of their fundamental groups.
	\end{example}

\section*{Acknowledgements}
The author wishes to express his thanks to Prof. Dr. Pavel Zalesskii for many discussions and advices.

\end{document}